\documentclass[12pt]{article}

\usepackage{a4wide,amsmath,amssymb,amsthm}
\usepackage{amsfonts}
\usepackage{fix-cm}
\usepackage{mathtools}
\usepackage{extarrows}
\usepackage{enumerate}
\usepackage{nicefrac}
\usepackage{afterpage}

\usepackage[utf8]{inputenc}
\usepackage[portuguese,main=english]{babel}
\usepackage[pdftex]{graphicx}	
\usepackage[authoryear]{natbib}
\usepackage[affil-it]{authblk}

\usepackage{geometry}
\usepackage[]{booktabs}
\usepackage{multirow,multicol}
\usepackage[flushleft]{threeparttable}
\usepackage{lscape}
\usepackage{longtable}

\usepackage{complexity}
\usepackage{marginnote}
\usepackage{textcomp}

\usepackage{url}
\usepackage{color}
\usepackage[x11names]{xcolor}
\colorlet{linkequation}{blue}
\usepackage[colorlinks]{hyperref}

\newcommand*{\SavedEqref}{}
\let\SavedEqref\eqref
\renewcommand*{\eqref}[1]{
\begingroup
\hypersetup{
    linkcolor=linkequation,
    linkbordercolor=linkequation,
}
\SavedEqref{#1}%
\endgroup
}

\hypersetup{
    citecolor=SpringGreen4,
    citebordercolor=SpringGreen4,
}

\allowdisplaybreaks
\usepackage{setspace}
\newtheorem{theorem}{Theorem}
\newtheorem{proposition}[theorem]{Proposition}

\theoremstyle{definition}

\theoremstyle{remark}

\newcommand{\defeq}{\stackrel{\text{def}}{=}}
\title{\Large 
On a class of strong valid inequalities for the connected matching polytope}
\author{Phillippe Samer}
\affil{\small 
{\tt $\lbrace$samer@uib.no$\rbrace$ }
}
\affil{\small 
Universitetet i Bergen, 
Institutt for Informatikk, 
Postboks 7800, 5020, Bergen, Norway
}

\date{September 23, 2023}

\begin{document}

\maketitle

\begin{abstract}
We identify a family of $O(|E(G)|^2)$ nontrivial facets of the connected matching polytope of a graph $G$, that is, the convex hull of incidence vectors of matchings in $G$ whose covered vertices induce a connected subgraph.
Accompanying software to further inspect the polytope of an input graph is available.
\end{abstract}

%%%%%%%%%%%%%%%%%%%%%%%%%%%%%%%%%%%%%%%%%%%%%%%%%%%%%%%%%%%%%%%%%%%%%%%%%%%%%%%%
%%%%%%%%%%%%%%%%%%%%%%%%%%%%%%%%%%%%%%%%%%%%%%%%%%%%%%%%%%%%%%%%%%%%%%%%%%%%%%%%
%%%%%%%%%%%%%%%%%%%%%%%%%%%%%%%%%%%%%%%%%%%%%%%%%%%%%%%%%%%%%%%%%%%%%%%%%%%%%%%%

\section{Introduction}
\label{sec:intro}

Our goal with this paper is to bring attention to an interesting polytope, and contribute towards improved algorithms for a recent combinatorial optimization problem defined over it.
%%% TO DO: highlight change in red
%A \textsf{P}-matching in a graph $G$ consists of a matching $M$ such that the subgraph induced by vertices covered by $M$ has property \textsf{P}, \textit{e.g.} being connected.
A \textsf{P}-matching in a graph $G$ is a matching $M$ such that the subgraph induced by vertices covered by $M$ has property \textsf{P}, \textit{e.g.} being connected.
In particular, while finding a maximum cardinality \textit{connected matching} is a well-solved problem, the edge-weighted counterpart is NP-hard even in very restricted graph classes~\citep{gg2023tcs}.
We initiate the systematic study of the  polytope $\mathfrak{C}(G)$ of connected matchings in a general graph $G$, introducing a relevant class of facet-defining inequalities.

Early examples of studies on \textsf{P}-matching problems include
\cite{StockmeyerVazirani1982} on induced matchings,
\cite{golumbic2001uniquely} on uniquely restricted matchings, and
\cite{Goddard2005disconnected} contemplating acyclic, connected and disconnected matchings.
While different \textsf{P}-matching problems receive increased attention in the recent literature, we highlight \cite{gg2022latin} and  \cite{gg2023tcs} on weighted connected matchings, who were able to determine several fine-grained complexity results. In particular, they show that it is NP-hard to find a maximum weight connected matching even in bounded-diameter bipartite and planar bipartite graphs.

The main argument for our contribution is to bring the machinery of polyhedral studies and mixed-integer programming (MIP) to bear on the investigation of weighted connected matchings in general graphs.
%Moreover, the combinatorial analysis of the facial structure of polytope $\mathfrak{C}$ is interesting in its own right.
In light of decades' worth of progress on matching theory and on the effective use of strong valid inequalities in MIP solvers, the combinatorial analysis of the facial structure of polytope $\mathfrak{C}(G)$ is a natural methodology.
On that perspective, the key idea we present next is a powerful ingredient in that direction.

We remark that all polyhedra in this work are rational, and that we do not leave the unit hypercube.
Nearly all terminology and notation we use are standard in graph theory and polyhedral combinatorics. 
The following might be worth mentioning.
We write $[k] \defeq \left\lbrace 1, \ldots, k \right\rbrace$.
Given a graph $G$, we denote its \textit{line graph} by $L(G)$, and define the \textit{distance between two edges} in $G$ as
$d_L: E(G) \times E(G) \rightarrow \mathbb{Z}_+$ so that $d_L(e_1, e_2)$ equals the number of edges in a shortest path between $e_1$ and $e_2$ in $L(G)$.
Given a subset of edges $S \subseteq E(G)$, we denote by $\chi^S$ its \textit{incidence} (or \textit{characteristic}) \textit{vector} in space $\mathbb{Q}^{|E(G)|}$, with $\chi^S_e = 1$ for each $e \in S$, and $\chi^S_e = 0$ otherwise.

First, it is convenient to clear implied equations out of systems of inequalities in studies of the connected matching polytope.
Using the standard argument that the unit vectors
$\left\lbrace \chi^{\lbrace e \rbrace}: e \in E(G) \right\rbrace$
and $\mathbf{0} \in \mathbb{Q}^{|E(G)|}$ are affinely independent and induce trivial incidence vectors of connected matchings in $G$, we have the following result.

\begin{proposition}
\label{thm:dim}
The connected matching polytope $\mathfrak{C}(G)$ of an arbitrary graph $G$ is full-dimensional.
\end{proposition}

Our main result is the following.
\begin{theorem}
\label{thm:facets}
Let $G$ be a connected graph and 
$\left\lbrace e^\prime, e^{\prime \prime} \right\rbrace \subset E(G)$
%%% TO DO: highlight change in red.
%induce a disconnected matching. Denote by 
be a disconnected matching. Denote by 
$\Lambda \defeq \Lambda(e^\prime, e^{\prime \prime}) = \left\lbrace f \in E(G): d_L(f,e^\prime) = d_L(f, e^{\prime \prime}) = 2 \right\rbrace$ the corresponding set of edges at two hops from both $e^{\prime}$ and $e^{\prime \prime}$.
%NB! this implies f \cap e^\prime = f \cap e^{\prime \prime} = \varnothing
%%% TO DO: highlight change in red.
%If $\Lambda \neq \varnothing$, the inequality
Suppose, further, that there is no connected matching including both $e^{\prime}$ and $e^{\prime \prime}$ in $G-\Lambda$ (\textit{i.e.} the subgraph of $G$ without edges in $\Lambda$). Then, the inequality
\vspace{-0.25cm}
%If $G-\Lambda$ (\textit{i.e.} the subgraph of $G$ without edges in $\Lambda$) is disconnected, the inequality
\begin{equation}
\label{eq:vi}
x_{e^\prime} + x_{e^{\prime \prime}} - \sum_{f \in \Lambda} x_f \leq 1
\end{equation}
is valid for $\mathfrak{C}(G)$.
Moreover, $\eqref{eq:vi}$ defines a facet when $\Lambda$ induces a clique in $L(G)$ and the subgraph induced by $\left\lbrace e^\prime, e^{\prime \prime}, f \right\rbrace$ is 2-connected for each $f \in \Lambda$.
\end{theorem}

The proof is saved for the next section.
Let us first illustrate how important it may be to consider this small set of facets (note that we have at most one inequality for each pair of edges).

There are many options for modelling induced connectivity.
Progress in mathematical programming computation of structures like maximum-weight connected subgraphs and Steiner trees build on vertex choosing binary variables $y \in \left\lbrace 0,1 \right\rbrace^{|V(G)|}$ and \textit{minimal separator inequalities} (MSI):
$
y_a + y_b -\sum_{u \in \mathcal{C}} y_u \leq 1
$
for each pair of non-adjacent vertices $a$ and $b$, and each $(a,b)$-separator $\mathcal{C} \subseteq V \backslash \left\lbrace a,b \right\rbrace$, \textit{i.e.} there are no paths connecting $a$ to $b$ if we remove $\mathcal{C}$ from $G$.
See \cite{wang2017imposing} for a thorough polyhedral analysis, and \cite{fischetti2017MPC} for supporting experimental results of an award-winning solver for Steiner tree problems.

In an attempt to build on those results, and impose induced connectivity on a system of inequalities formulating connected matchings  while using only natural design variables $x \in \left\lbrace 0,1 \right\rbrace^{|E(G)|}$, as opposed to working with extended formulations, one may use the fact that vertex $u$ belongs to the subgraph induced by matching $M$ if and only if there is exactly one edge in $M$ incident to $u$. That is, projecting MSI onto the space of $x$ variables using $y_u \defeq \sum_{e \in \delta(u)} x_e$, we derive a MIP formulation to find maximum weight connected matchings using MSI.
We are currently pursuing that endeavour and crafting a branch-and-cut algorithm for weighted connected matchings.
In the meantime, we inspected the convex hull $\mathfrak{C}(G)$ for several examples using \textsf{polymake}~\citep{polymake2000,polymake2017}.
Remarkably, we discovered fine examples where a single inequality in~$\eqref{eq:vi}$ dominates several MSI.

\begin{figure}[t]
\begin{center}
    \includegraphics[scale=0.6]{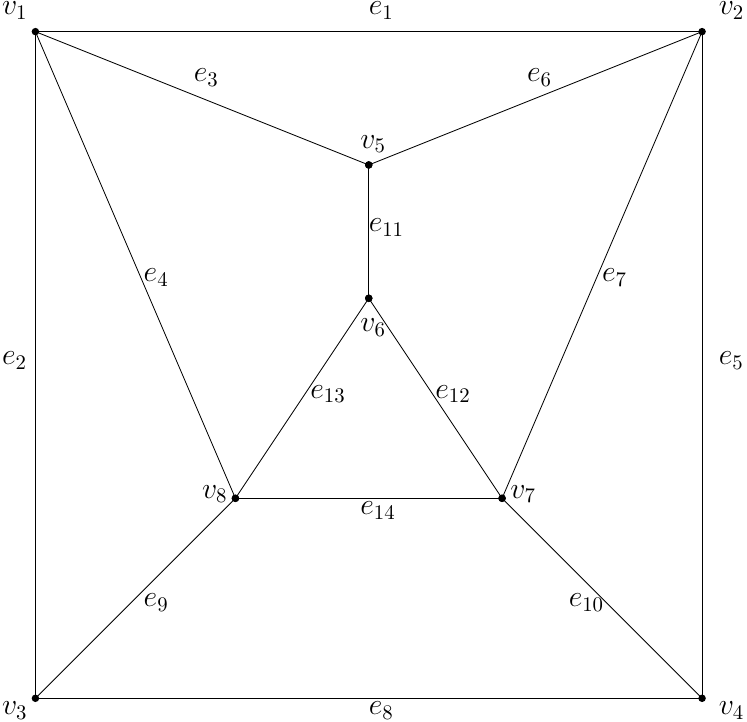}
\end{center}
\caption{Skeleton graph of Johnson Solid 26.}
\label{fig:JS26}
\end{figure}

For instance, taking $G_{J26}$ to be the skeleton graph of Johnson Solid 26, depicted in Figure~\ref{fig:JS26}, and studying the minimal inequality description of $\mathfrak{C}(G_{J26})$, we detect 14 trivial facets from non-negativity bounds,
8 blossom inequalities on handles $H_i=V(G)\backslash \left\lbrace v_i \right\rbrace$, some blossom inequalities on triangles, and 5 facets defined by our inequalities in~$\eqref{eq:vi}$.
Among the latter, we find
\begin{equation}
\label{eq:example_J26}
x_5 + x_{13} -x_2 -x_3 \leq 1,
\end{equation}
which may be interpreted as the lifting of 4 different MSI corresponding to $\mathcal{C} \defeq \left\lbrace v_1, v_3, v_5, v_7 \right\rbrace$ as a minimal $(v_a, v_b)$-separator for
$(a,b) \in \left\lbrace
(2,6), (2,8), (4,6), (4,8)
\right\rbrace$.
In other words, a MIP formulation adding inequality~$\eqref{eq:example_J26}$ \textit{a priori} gives a tighter approximation of $\mathfrak{C}(G_{J26})$ than a formulation that depends solely on cutting planes from MSI projected onto $x$ space, which could generate dynamically the cuts
$
\left\lbrace y_{v_a} +y_{v_b} -\sum_{z \in \mathcal{C}} y_{z} \leq 1 \right\rbrace
$.
For example, when $(a,b) = (2,6)$, the projected inequality is
\begin{footnotesize}
\begin{alignat}{2}
 (x_1 + x_5 + x_6 + x_7)
+(x_{11}+x_{12}+x_{13}) & \nonumber \\
-(x_1+x_2+x_3+x_4)
-(x_2+x_8+x_9) & \nonumber \\
-(x_3+x_6+x_{11})
-(x_7+x_{10}+x_{12}+x_{14})
& \leq 1, \nonumber
\end{alignat}
\end{footnotesize}
%%% TO DO: REMOVE THIS ON DIFFERENT TYPESETTING
\hspace{-0.13cm}that is,
\begin{footnotesize}
$
 x_5 +x_{13}
-2x_2
-2x_3
-x_4
-x_8
-x_9
-x_{10}
-x_{14}
\leq 1,
$
\end{footnotesize}
which is immediately seem to be dominated by~$\eqref{eq:example_J26}$.

The productive exercise of plugging different input graphs $G$ and inspecting the resulting polytope $\mathfrak{C}(G)$ is made free and open to anyone interested.
The code developed to find a $\mathcal{V}$-description of the polytope corresponding to an input graph, and then produce an $\mathcal{H}$-description with \textsf{polymake}, is available at 
\begin{small}
\url{https://github.com/phillippesamer/wcm-branch-and-cut/tree/main/polyhedra},
\end{small}
as is the forthcoming branch-and-cut algorithm.

%%% TO DO: HIGHLIGHT CHANGE IN RED
%Before we resume to the main proof, we remark that it may be straightforward in some particular cases to show that $\eqref{eq:vi}$ induces a facet of the smaller polytope $\mathfrak{C}(G[S])$, with $S = \left\lbrace e^\prime, e^{\prime \prime}, f \right\rbrace$ and $f \in \Lambda$.
Before we resume to the main proof, we note that it may be straightforward in some particular cases to show that $\eqref{eq:vi}$ induces a facet of the smaller polytope $\mathfrak{C}(G[S])$, with $S = \left\lbrace e^\prime, e^{\prime \prime}, f \right\rbrace$ and $f \in \Lambda$.
Nonetheless, this would only be interesting if coupled with a lifting result showing how to derive a facet of our target polytope $\mathfrak{C}(G)$ from the smaller dimensional one.
We skip that altogether and give next a direct proof of the general result.

It is also important to remark that there are simple sufficient conditions to verify the theorem hypothesis that no connected matching including both $e^\prime$ and $e^{\prime \prime}$ exists in $G-\Lambda$.
For instance, it suffices that $G-\Lambda$ does not contain a path joining the two edges.
The latter is not, however, a necessary condition, as one can see using again the example underlying inequality~$\eqref{eq:example_J26}$ above.
We leave the problem of finding a strict condition for future work, and conjecture that it may be verified efficiently by working on the subgraph obtained from $G$ after removing $\Lambda$, as well as those edges adjacent to either $e^\prime$ or  $e^{\prime \prime}$, and using the fact that finding a maximum cardinality connected matching can be done in polytime.

%%%%%%%%%%%%%%%%%%%%%%%%%%%%%%%%%%%%%%%%%%%%%%%%%%%%%%%%%%%%%%%%%%%%%%%%%%%%%%%%
%%%%%%%%%%%%%%%%%%%%%%%%%%%%%%%%%%%%%%%%%%%%%%%%%%%%%%%%%%%%%%%%%%%%%%%%%%%%%%%%
%%%%%%%%%%%%%%%%%%%%%%%%%%%%%%%%%%%%%%%%%%%%%%%%%%%%%%%%%%%%%%%%%%%%%%%%%%%%%%%%

\section{Proof of Theorem~\ref{thm:facets}}

\paragraph{}\textit{(i)}
%%% TO DO: highlight changes in red
%The validity of inequality $\eqref{eq:vi}$ follows from simple combinatorial reasoning.
%Let $\chi^{S}$ be the incidence vector of an arbitrary connected matching $S$.
%If $e^\prime \not\in S$ or $e^{\prime \prime} \not\in S$, the left-hand side of the expression in $\eqref{eq:vi}$ is at most 1.
%On the other hand, the fact that $\left\lbrace e^\prime, e^{\prime \prime} \right\rbrace$ induces a disconnected matching implies that,   if both $e^\prime \in S$ and $e^{\prime \prime} \in S$, there must be an additional edge $f$ in $S$ to connect their endpoints. 
%
%From $\left\lbrace e^\prime, e^{\prime \prime}, f \right\rbrace$ being a matching, we have that $d_L(f,e^\prime) \geq 2$ and $d_L(f, e^{\prime \prime}) \geq 2$.
%For the corresponding vertices in each of $e^\prime$, $e^{\prime \prime}$ and $f$ to induce a connected subgraph, we verify $d_L(f,e^\prime) \leq 2$ and $d_L(f, e^{\prime \prime}) \leq 2$.
%Hence $f \in \Lambda$, and again the left-hand side of the inequality is bound below $1$ at $\chi^{S}$.
%Since $\chi^{S}$ is taken as an arbitrary vertex of the polytope, the inequality is valid at all points in $\mathfrak{C}(G)$ by convexity.
%
The validity of inequality $\eqref{eq:vi}$ follows from simple combinatorial reasoning.
Let $\chi^{S}$ be the incidence vector of an arbitrary connected matching $S$.
If $e^\prime \not\in S$ or $e^{\prime \prime} \not\in S$, the left-hand side of the expression in $\eqref{eq:vi}$ is at most 1.
On the other hand,
if both $e^\prime \in S$ and $e^{\prime \prime} \in S$, the hypothesis that there is no connected matching including both $e^\prime$ and $e^{\prime \prime}$ in $G-\Lambda$ implies that there must be an additional edge $f \in \Lambda$ in $S$, and again the left-hand side of the inequality is bound below $1$ at $\chi^{S}$.
Since $\chi^{S}$ is taken as an arbitrary vertex of the polytope, the inequality is valid at all points in $\mathfrak{C}(G)$ by convexity.

\paragraph{}\textit{(ii)}
For the facet proof, let
%Next, we use the facet-proving technique based on the maximality of a face \citep[Theorem I.4.3.6]{NemhauserWolsey1999}. Namely, let
$F \defeq \left\lbrace x \in \mathfrak{C}(G): x_{e^\prime} + x_{e^{\prime \prime}} - \sum_{f \in \Lambda} x_f = 1 \right\rbrace = \left\lbrace x \in \mathfrak{C}(G): \pi x = \pi_0\right\rbrace$
denote the face corresponding to inequality $\eqref{eq:vi}$; vector $(\pi, \pi_0)$ is just shorthand notation here.
By the full-dimensionality observation in Proposition~\ref{thm:dim}, there is no equation satisfied by all points in the polytope.
%%% TO DO: highlight change in red
%It thus suffice to show that if
It thus suffices to show that if
$F \subseteq \overline F \defeq \left\lbrace x \in \mathfrak{C}(G): \lambda x = \lambda_0 \right\rbrace$,
the defining inequalities $\pi x \leq \pi_0$ and $\lambda x \leq \lambda_0$ are actually equivalent (\textit{i.e.} there exists $\rho > 0$ such that $\lambda = \rho \pi$ and $\lambda_0 = \rho \pi_0$), and hence the strict inclusion cannot hold \citep[Theorem I.4.3.6]{NemhauserWolsey1999}.

Set $\rho \defeq \lambda_0$. To prove the equivalence of the nonzero coefficients, let us denote by 
$x^1, x^2, x^3$ 
the incidence vectors of the connected matchings consisting of
$\left\lbrace e^\prime \right\rbrace$,
$\left\lbrace e^{\prime \prime} \right\rbrace$, and
$\left\lbrace e^\prime, e^{\prime \prime}, f_1 \right\rbrace$, respectively, with $f_1 \in \Lambda$ arbitrary. 
Note that each of these points belong to $F$, and give a single nonzero coefficient when evaluating $\lambda x = \lambda_0$:
\begin{enumerate}[i.]
\item From 
$x^{1} \defeq \left\lbrace x^{1}_{e^\prime} = 1, x^{1}_{*} = 0 \text{ for } * \in E\backslash \left\lbrace e^\prime \right\rbrace \right\rbrace 
\in F$
we get
$\lambda x^{1} = \lambda_{e^{\prime}}$.
Together with $\lambda x^{1} = \lambda_0$, we obtain
$\lambda_{e^{\prime}} = \lambda_0 \defeq \rho = \rho \cdot 1 = \rho \cdot \pi_{e^{\prime}}$.
%%%%%%%%%%
\item From 
$x^{2} \defeq \left\lbrace x^{2}_{e^{\prime \prime}} = 1, x^{2}_{*} = 0 \text{ for } * \in E\backslash \left\lbrace e^{\prime \prime} \right\rbrace \right\rbrace 
\in F$
we get
$\lambda x^{2} = \lambda_{e^{\prime \prime}}$.
Together with $\lambda x^{2} = \lambda_0$, we obtain
$\lambda_{e^{\prime \prime}} = \lambda_0 \defeq \rho = \rho \cdot 1 = \rho \cdot \pi_{e^{\prime \prime}}$.
%%%%%%%%%%
\item From
$x^{3} \defeq \left\lbrace x^{3}_{e^\prime} = x^{3}_{e^{\prime \prime}} = x^{3}_{f_1} = 1, x^{3}_{*} = 0 \text{ for } * \in E\backslash \left\lbrace e^\prime, e^{\prime \prime}, f_1 \right\rbrace \right\rbrace 
\in F$
we get
$\lambda x^{3} = \lambda_{e^{\prime}} + \lambda_{e^{\prime \prime}} + \lambda_{f_1}$.
Together with $\lambda x^{3} = \lambda_0$, we have
$\lambda_{f_1} = \lambda_0 -\lambda_{e^{\prime}}  -\lambda_{e^{\prime \prime}}$.
Now, using the coefficients determined in the previous two items, we find that 
$\lambda_{f_1} = \lambda_0 -\lambda_0 -\lambda_0 = \lambda_0 \cdot (-1) \defeq \rho \cdot (-1) = \rho \cdot \pi_{f_1}$.
\end{enumerate}
Since the choice of $f_1 \in \Lambda$ is without loss of generality, we repeat the argument in the last item above using the incidence vectors of each of the connected matchings in
$
\left\lbrace
\left\lbrace e^\prime, e^{\prime \prime}, f \right\rbrace : f \in \Lambda
\right\rbrace
$,
to actually determine that
$\lambda_{f} = \rho \cdot \pi_f$ for each $f \in \Lambda$.

That establishes the equivalence of nonzero coefficients. It remains to show that all other coefficients of $\lambda$ are null.
All we need for that are the following two remarks.

%%%%%%%%%%%%%%%%%%%%%%%%%%%%%%%%%%%%%%%%%%%%%%
%%%%%%%%%%%%%%%%%%%%%%%%%%%%%%%%%%%%%%%%%%%%%%
\begin{description}

%%%%%%%%%%%%%%%%%%%%%%%%%%%%%%%%%%%%%%%%%%%%%%
\item[Fact 1] Let $\chi^M \in F$ be the incidence vector of matching $M$.
By the hypothesis that $\Lambda$ induces a clique in $L(G)$, we have that $M$ contains at most one edge in~$\Lambda$.
That enables us to judiciously order the edges in $M = \left\lbrace e_1, \ldots, e_m \right\rbrace$ so that
(i) the edges in $\overline{\Lambda} \defeq \left\lbrace e^\prime, e^{\prime \prime}\right\rbrace \cup \Lambda$ (either one or exactly three) appear first, and
(ii) each edge added after we are done with $\overline{\Lambda}$ yields an additional point in our face $F$.
Ultimately, we may apply the following simple observation a number of times:
\begin{alignat}{2}
\chi^{\left\lbrace e_1, \ldots, e_k \right\rbrace} \in F \subseteq \overline{F}
\implies &
\lambda \cdot \chi^{\left\lbrace e_1, \ldots, e_k \right\rbrace} =
\overbrace{\sum_{i=1}^{k-1} \lambda_{e_i} \cdot 1}^{= \lambda_0} \ 
+ \lambda_{e_k} \cdot 1 \ 
+ \overbrace{\sum_{\substack{e \in E(G): \\e\not\in\left\lbrace e_1, \ldots, e_k \right\rbrace}} \lambda_e \cdot 0}^{= 0} \ 
= \ \lambda_0 \nonumber \\
\implies &
\lambda_{e_k} = 0 \label{eq:augmenting_lambda}
\end{alignat}

%%%%%%%%%%%%%%%%%%%%%%%%%%%%%%%%%%%%%%%%%%%%%%
%%% TO DO: highlight change in red
%\item[Fact 2] For all $e \in E(G)$, there exists a connect matching $M$ including $e$ such that $\chi^M \in F$.
\item[Fact 2] For all $e \in E(G)$, there exists a connected matching $M$ including $e$ such that $\chi^M \in F$.
Let $S$ be the set of endpoints of edges in $\overline{\Lambda} = \left\lbrace e^\prime, e^{\prime \prime}\right\rbrace \cup \Lambda$, and consider the following two cases.
\begin{enumerate}[(i)]
\item Suppose that $e \in E(G[S])$.
Taking $M$ in
$
\left\lbrace
\left\lbrace e^\prime, e^{\prime \prime}, f \right\rbrace : f \in \Lambda
\right\rbrace
$, the claim follows immediately for $e \in \overline{\Lambda}$, by definition of $F$.
Otherwise, if $e = \left\lbrace u,v \right\rbrace$ is an induced edge in
$E(G[S])\backslash \overline{\Lambda}$, it joins a vertex from either $e^\prime$ or $e^{\prime \prime}$ to a vertex of some $f \in \Lambda$.
Suppose without loss of generality that $u \in e^\prime$ and $v \in f$.
Since the subgraph induced by $\left\lbrace e^\prime, e^{\prime \prime}, f \right\rbrace$ is 2-connected by the hypothesis in the theorem, the matching $M \defeq \left\lbrace e, e^{\prime \prime} \right\rbrace$ is connected, and $\chi^M \in F$.

\item If $\hat{e} \not\in E(G[S])$, we use the assumption that $G$ is connected and take a simple path $P = (e_1, \ldots, e_n)$ on $n$ edges beginning at $e_1 = \hat{e}$, and minimal with respect to including an edge in $E(G[S])$, that is, $e_n$ is the only edge in $P$ which is also in the subgraph induced by vertices covered by $\left\lbrace e^\prime, e^{\prime \prime}, \hat{f} \right\rbrace$, for some $\hat{f} \in \Lambda$.
%%% TO DO: highlight change in red: comma after Now
Now, we let $\hat{M}$ be the matching alternating edges along $P$, \textit{requiring} $\hat{e} \in \hat{M}$, and reason on the parity of the path length.

Suppose first that $n$ is odd, so that $e_n$ is also in $\hat{M}$.
If $e_n$ is either $e^\prime$ or $e^{\prime \prime}$, we are done.
If $e_n = \hat{f}$, we set $M \defeq \hat{M} \uplus \left\lbrace e^\prime, e^{\prime \prime}\right\rbrace$, and we are done.
Otherwise, $e_n = \left\lbrace u,v \right\rbrace$ as in case (i): without loss of generality, let $e_n$ join $u \in e^\prime$ and $v \in \hat{f}$.
Using again the hypothesis that the subgraph induced by $\left\lbrace e^\prime, e^{\prime \prime}, \hat{f} \right\rbrace$ is 2-connected, choose $M \defeq \hat{M} \uplus \left\lbrace e^{\prime \prime} \right\rbrace$, and we are done.

Suppose now that $n$ is even, and thus $e_n \not\in \hat{M}$.
By the 2-connectivity hypothesis, we may augment $\hat{M}$ to get a connected matching with exactly one of $e^\prime$ or $e^{\prime \prime}$, whose characteristic vector is thus in the face $F$. 
If $e_{n-1}$ is incident to a vertex of $e^\prime$, we set $M \defeq \hat{M} \uplus \left\lbrace e^{\prime \prime}, e^\dag \right\rbrace$, where $e^\dag$ joins the endpoint of $e^\prime$ not covered by~$\hat{M}$ to a vertex of $\hat{f}$.
The analogous argument holds when $e_{n-1}$ is incident to a vertex of $e^{\prime \prime}$.
Finally, if $e_{n-1}$ is incident to a vertex of $\hat{f}$, we simply choose $M \defeq \hat{M} \uplus \left\lbrace e^\prime \right\rbrace$.
\end{enumerate}

\end{description}
%%%%%%%%%%%%%%%%%%%%%%%%%%%%%%%%%%%%%%%%%%%%%%
%%%%%%%%%%%%%%%%%%%%%%%%%%%%%%%%%%%%%%%%%%%%%%

To complete the proof, we consider each edge $e_k \not\in \overline{\Lambda}$.
By Fact~2, we obtain a connected matching $M$ including $e_k$ such that $\chi^M \in F$.
Using~$\eqref{eq:augmenting_lambda}$ in Fact~1, we conclude that $\lambda_{e_k} = 0$.
Hence $(\lambda, \lambda_0) = (\rho \pi, \rho \pi_0)$, and $F$ determines a facet of $\mathfrak{C}(G)$.
\qed

%%%%%%%%%%%%%%%%%%%%%%%%%%%%%%%%%%%%%%%%%%%%%%%%%%%%%%%%%%%%%%%%%%%%%%%%%%%%%%%%
%%%%%%%%%%%%%%%%%%%%%%%%%%%%%%%%%%%%%%%%%%%%%%%%%%%%%%%%%%%%%%%%%%%%%%%%%%%%%%%%

\paragraph{Acknowledgement} 
The author is grateful to the support by the Research Council of Norway through the research project 249994 CLASSIS. This work is dedicated to the sweet memory of our department administration member Ingrid Kyllingmark.

%\bibliographystyle{plainnat-no-doi}
%\bibliography{references}

\begin{thebibliography}{10}
\providecommand{\natexlab}[1]{#1}
\providecommand{\url}[1]{\texttt{#1}}
\expandafter\ifx\csname urlstyle\endcsname\relax
  \providecommand{\doi}[1]{doi: #1}\else
  \providecommand{\doi}{doi: \begingroup \urlstyle{rm}\Url}\fi

\bibitem[Assarf et~al.(2017)Assarf, Gawrilow, Herr, Joswig, Lorenz, Paffenholz,
  and Rehn]{polymake2017}
B.~Assarf, E.~Gawrilow, K.~Herr, M.~Joswig, B.~Lorenz, A.~Paffenholz, and
  T.~Rehn.
\newblock Computing convex hulls and counting integer points with polymake.
\newblock \emph{Mathematical Programming Computation}, 9:\penalty0 1--38, 2017.
\newblock URL \url{https://doi.org/10.1007/s12532-016-0104-z}.

\bibitem[Fischetti et~al.(2017)Fischetti, Leitner, Ljubi{\'c}, Luipersbeck,
  Monaci, Resch, Salvagnin, and Sinnl]{fischetti2017MPC}
M.~Fischetti, M.~Leitner, I.~Ljubi{\'c}, M.~Luipersbeck, M.~Monaci, M.~Resch,
  D.~Salvagnin, and M.~Sinnl.
\newblock Thinning out steiner trees: a node-based model for uniform edge
  costs.
\newblock \emph{Mathematical Programming Computation}, 9\penalty0 (2):\penalty0
  203--229, 2017.
\newblock URL \url{https://doi.org/10.1007/s12532-016-0111-0}.

\bibitem[Gawrilow and Joswig(2000)]{polymake2000}
E.~Gawrilow and M.~Joswig.
\newblock Polymake: a framework for analyzing convex polytopes.
\newblock In G.~Kalai and G.~M. Ziegler, editors,
  \emph{Polytopes—combinatorics and computation}, pages 43--73. Springer,
  2000.
\newblock URL \url{https://doi.org/10.1007/978-3-0348-8438-9_2}.

\bibitem[Goddard et~al.(2005)Goddard, Hedetniemi, Hedetniemi, and
  Laskar]{Goddard2005disconnected}
W.~Goddard, S.~M. Hedetniemi, S.~T. Hedetniemi, and R.~Laskar.
\newblock Generalized subgraph-restricted matchings in graphs.
\newblock \emph{Discrete Mathematics}, 293\penalty0 (1):\penalty0 129--138,
  2005.
\newblock ISSN 0012-365X.
\newblock URL \url{https://doi.org/10.1016/j.disc.2004.08.027}.
\newblock 19th British Combinatorial Conference.

\bibitem[Golumbic et~al.(2001)Golumbic, Hirst, and
  Lewenstein]{golumbic2001uniquely}
M.~C. Golumbic, T.~Hirst, and M.~Lewenstein.
\newblock Uniquely restricted matchings.
\newblock \emph{Algorithmica}, 31:\penalty0 139--154, 2001.
\newblock URL \url{https://doi.org/10.1007/s00453-001-0004-z}.

\bibitem[Gomes et~al.(2022)Gomes, Masquio, Pinto, Santos, and
  Szwarcfiter]{gg2022latin}
G.~C.~M. Gomes, B.~P. Masquio, P.~E.~D. Pinto, V.~F.~d. Santos, and J.~L.
  Szwarcfiter.
\newblock Weighted connected matchings.
\newblock In A.~Casta{\~{n}}eda and F.~Rodr{\'i}guez-Henr{\'i}quez, editors,
  \emph{LATIN 2022: Theoretical Informatics}, pages 54--70, Cham, 2022.
  Springer International Publishing.
\newblock ISBN 978-3-031-20624-5.
\newblock URL \url{https://doi.org/10.1007/978-3-031-20624-5_4}.

\bibitem[Gomes et~al.(2023)Gomes, Masquio, Pinto, {dos Santos}, and
  Szwarcfiter]{gg2023tcs}
G.~C. Gomes, B.~P. Masquio, P.~E. Pinto, V.~F. {dos Santos}, and J.~L.
  Szwarcfiter.
\newblock Disconnected matchings.
\newblock \emph{Theoretical Computer Science}, 956:\penalty0 113821, 2023.
\newblock ISSN 0304-3975.
\newblock URL \url{https://doi.org/10.1016/j.tcs.2023.113821}.

\bibitem[Nemhauser and Wolsey(1999)]{NemhauserWolsey1999}
G.~L. Nemhauser and L.~A. Wolsey.
\newblock \emph{Integer and combinatorial optimization}, volume~55 of
  \emph{Wiley-Interscience series in discrete mathematics and optimization}.
\newblock John Wiley \& Sons, Inc, 1999.
\newblock URL \url{https://doi.org/10.1002/9781118627372}.

\bibitem[Stockmeyer and Vazirani(1982)]{StockmeyerVazirani1982}
L.~J. Stockmeyer and V.~V. Vazirani.
\newblock {NP-completeness of some generalizations of the maximum matching
  problem}.
\newblock \emph{Information Processing Letters}, 15\penalty0 (1):\penalty0
  14--19, 1982.
\newblock URL \url{https://doi.org/10.1016/0020-0190(82)90077-1}.

\bibitem[Wang et~al.(2017)Wang, Buchanan, and Butenko]{wang2017imposing}
Y.~Wang, A.~Buchanan, and S.~Butenko.
\newblock On imposing connectivity constraints in integer programs.
\newblock \emph{Mathematical Programming}, 166:\penalty0 241--271, 2017.
\newblock URL \url{https://doi.org/10.1007/s10107-017-1117-8}.

\end{thebibliography}

\end{document}